\documentclass[epic,eepic,11pt]{amsart}

\usepackage{amsfonts,mathrsfs}
\usepackage[mathscr]{eucal}

\usepackage{amssymb}
\usepackage{amscd}
\usepackage{fancyhdr}

\usepackage{euler,eucal}

\DeclareMathAlphabet{\mathbf}{T1}{ppl}{bx}{n}
\DeclareMathAlphabet{\mathrm}{T1}{ppl}{m}{n}

\numberwithin{equation}{section}

\newcommand\note[1]%
{$^\dagger$\marginpar{\footnotesize{$^\dagger${#1}}}}

\def\({\left(}
\def\){\right)}
\def\<{\left<}
\def\>{\right>}

\newtheorem{theorem}{Theorem}[section]
\newtheorem{proposition}[theorem]{Proposition}
\newtheorem{lemma}[theorem]{Lemma}
\newtheorem{definition}[theorem]{Definition}

\theoremstyle{definition}
\newtheorem{example}[theorem]{Example}

\newcommand\bb[1]{{\text{\bf#1}}}

\newcommand\Z{\bb{Z}}

\newcommand\R{\mathbb{R}}

\renewcommand\H{\bb{H}}

\newcommand\funclim[1]{\operatorname*{\mathrm{#1}}}

\renewcommand\lim{\funclim{lim}}

\newcommand\sur{\mathrel{\to\kern-1.8ex\to}}
\newcommand\iso{\mathrel{\hookrightarrow\kern-1.8ex\to}}

\newcommand\longhookrightarrow{\lhook\joinrel\longrightarrow}

\newcommand\longsur{\mathrel{\longrightarrow\kern-1.8ex\to}}
\newcommand\longiso{\mathrel{\longhookrightarrow\kern-1.8ex\to}}

\begin{document}

\bibliographystyle{amsalpha}
\date{\today}

\title{Examples of non-K\"ahler Hamiltonian circle manifolds with the strong
Lefschetz property}

\author{Yi Lin}

\begin{abstract}
In this paper we construct six-dimensional compact non-K\"ahler
Hamiltonian circle manifolds which satisfy the strong Lefschetz
property themselves but nevertheless have a non-Lefschetz symplectic
quotient. This provides the first known counter examples to the
question whether the strong Lefschetz property descends to the
symplectic quotient. We also give examples of Hamiltonian strong
Lefschetz circle manifolds which have a non-Lefschetz fixed point
submanifold. In addition, we establish a sufficient and necessary
condition for a finitely presentable group to be the fundamental
group of a strong Lefschetz manifold. We then use it to show the
existence of Lefschetz four-manifolds with non-Lefschetz finite
covering spaces.
\end{abstract}

\maketitle

\section{Introduction}

   Brylinski defined in \cite{brylinski;differential-poisson} the
   notion of symplectic harmonic forms. He further conjectured
   that on a compact symplectic manifold every cohomology class has a
   harmonic representative and proved that this is the case for compact
K\"ahler
   manifolds and certain other examples.

    A symplectic manifold $(M, \omega)$ of dimension $2m$ is said to have
the strong
    Lefschetz property or equivalently to be a strong Lefschetz manifold if
    and only if for any $0 \leqslant k \leqslant m$, the Lefschetz type map
    \begin{equation}\label{lefschetz property}
    L^k_{[\omega]} : H^{m-k}(M) \rightarrow H^{m+k}(M),\,\,
    [\alpha] \rightarrow [\alpha \wedge \omega^k]
    \end{equation}
    is onto. Mathieu  \cite{mathieu;harmonic-symplectic} proved the
    remarkable theorem that Brylinski conjecture is true for
    a symplectic manifold $(M, \omega)$ if and only if it has
    the strong Lefschetz property. This result was strengthened
    by Merkulov  \cite{merkulov;formality-canonical-symplectic}
    and Guillemin \cite{guillemin;symplectic-hodge},
    who independently established the symplectic $d, \delta$-lemma for
compact
    symplectic manifolds with the strong Lefschetz property. As a consequence
    of the symplectic $d, \delta$-lemma, they showed that strong
     Lefschetz manifolds are formal in a certain sense.

    We obtained an equivariant version of the above results jointly
    with Sjamaar in \cite{LS}. In particular, it was proved in \cite{LS} that for a compact Hamiltonian
$G$-manifold with the strong Lefschetz property every cohomology
class has a canonical equivariant extension. In a
    subsequent paper \cite{Lin} the author extended the main results in
\cite{LS}
    to equivariant differential forms with generalized coefficients on
Hamiltonian
    manifolds with the strong Lefschetz property.

    Kaoru Ono and Reyer Sjamaar raised the question whether the strong
Lefschetz property descends
    to the symplectic quotient. Obviously, this question has an affirmative
answer in the category of
     equivairant K\"ahler geometry. It is then a very natural question to
ask whether
     this is still the case in the symplectic category.

   The main result of this paper is first known counter examples
   which show,  in contrast with the equivariant K\"ahler case, that
   the strong Lefschetz property does not survive
    symplectic reduction in general.   The difficulty constructing such
examples comes largely from
    the lack of general examples of non-K\"ahler  Hamiltonian symplectic
manifolds which have the strong
    Lefschetz property. Historically, a lot of examples of non-K\"ahler
symplectic manifolds have
  been constructed.  However, as the strong Lefschetz property is commonly
used as a tool to detect
  the existence of K\"ahler structure, not many known examples of
non-K\"ahler
  symplectic manifolds have the strong Lefschetz property.

   By Mathieu's theorem \cite{mathieu;harmonic-symplectic}
for a symplecitc manifold with the strong Lefschetz property the
symplectic harmonic groups always coincide with the de Rham
cohomology groups.  It is noteworthy that Dong
\cite{yan;hodge-structure-symplectic} showed that there exist
compact symplectic four-manifolds which admit a family $\omega_t$ of
symplectic forms such that the dimension
    of the third symplectic harmonic group varies.  Dong's construction
depends heavily on the following
    result in Gompf's path-breaking paper \cite{R. Gompf}.

\begin{theorem} \label{Gompf}\footnote{The first assertion of Theorem
\ref{Gompf} is contained in the
     statement of Theorem 4.1 of \cite{R. Gompf}; the second
     assertion follows from the discussion following the proof of
     Observation 7.4 in the same paper.}(\cite{R. Gompf}) Let $G$ be any
finitely presentable group. Then there is a closed, symplectic
$4$-manifold $(M, \omega)$ such that \begin{enumerate}
\item [(1)] $\pi _1(M)=G$, \item [(2)] The Lefschetz map
$L_{[\omega]}: H^1(M) \rightarrow H^3(M)$ is trivial.\end{enumerate}
\end{theorem}

    Our construction of counter examples is inspired by the above-mentioned
work of
    Dong  \cite{yan;hodge-structure-symplectic}, and by Karshon's
    example \cite{Yael}  of a Hamiltonian circle six-manifold with
    a non-log concave Duistermaat-Heckman function, which in turn is
    a  piece of a manifold constructed by Mcduff \cite{Mcduff}.
    Let us give a brief account of the main ideas of our construction
    here.   First, we show that any finitely presentable group $G$ with
   a certain structure can be realized as the fundamental group of a
    four-manifold $N$ which supports a family of symplectic forms $\omega_t,
t\in \R$, such that $(N, \omega_0)$
    does not have the strong Lefschetz property. Second, we
     prove that for such a manifold $N$ there exists a six-dimensional
compact
    Hamiltonian symplectic $S^1$-manifold $M$ which is fibred over $N$ with
the fibre $S^2$; furthermore,
    the symplectic quotient of $M$ taken at a certain value will be exactly
$N$ with the reduced
    form $\omega_{0}$.  As $G$ varies, we actually obtain infinitely many
topologically inequivalent six-dimensional compact Hamiltonian
    $S^1$-manifolds, each of which has the strong Lefschetz property itself
but nevertheless admits a non-Lefschetz
     symplectic quotient.  This also gives us new examples of  compact
non-K\"ahler Hamiltonian manifolds. (C.f.,  \cite{E. Lerman2} and
\cite{S. Tolman}.)

   The same ideas also allow us to construct Hamiltonian strong
Lefschetz manifolds which have non-Lefschetz fixed point
submanifold. For a compact Hamiltonian manifold, one interesting
question is what the relationship is between the symplectic harmonic
theory of the manifold itself and that of its fixed point
submanifold. For instance it remains an open question whether a
compact Hamiltonian circle manifold with isolated fixed points has
to satisfy the strong Lefschetz property. And one may further ask
whether the strong Lefschetz property for a Hamiltonian circle
manifold and its fixed point submanifolds will imply each other. Our
examples give a negative answer to the latter question.

As an aside, we give a sufficient and necessary condition for a
finitely presentable group to be the fundamental group of a compact
symplectic four manifold with the strong Lefshetz property. It
suggests that the fundamental groups of strong Lefschetz manifolds
and that of K\"ahler manifolds may have quite different behavior. In
fact, it enables us to construct examples of compact strong
Lefschetz manifolds which have non-Lefschetz finite covering spaces.

It is an important question to which extent the symplectic manifolds
are more general than  K\"ahler manifolds. The examples constructed
in this paper show clearly that the category of strong Lefschetz
manifolds with Hamiltonian circle actions is much larger than the
category of K\"ahler  manifolds with compatible Hamiltonian circle
actions.

    This paper is organized as follows.
    Section \ref{examples} modifies Dong Yan's methods
\cite{yan;hodge-structure-symplectic}
    to prove the existence of the symplectic four-manifolds with certain
properties we want.
     Section \ref{fundamental group} records a sufficient and
   necessary condition for a finitely presentable group $G$ to be
   the fundamental group of a compact  strong Lefschetz four manifold. As an
immediate application
   of this observation, Section \ref{fundamental group} also gives us
examples of strong
   Lefschetz manifolds with non-Lefschetz finite covering spaces. Section
\ref{construction} shows how to construct compact Hamiltonian strong
Lefschetz circle manifolds with a
   non-Lefschetz symplectic quotient. In addition, Section \ref{construction}
   also explains how to obtain examples of Hamiltonian Strong Lefschetz circle
manifolds
   with a non-Lefschetz fixed point submanifold.

    \,\,\,\,\,\,\,\,\,\,\,\,\,\,\,\,\,

   \textbf{Acknowledgement}:  I would like  to thank  Reyer Sjamaar for
posing to me the question
    if the strong Lefschetz property survives the symplectic reduction,
for pointing out to me
    Karshon's example \cite{Yael}, and for many stimulating discussions on strong Lefschetz property in general.

    I would like to thank Kenneth Brown for explaining to me some formulas
in group
    cohomology. And I would also like to thank the referees for many
    valuable advices.

   A special thank goes to Victor Guillemin, who advised me to
   take a close look at Dong Yan's examples in
\cite{yan;hodge-structure-symplectic}.
   The present paper would not have been possible without his prescient
advice.


\section{Symplectic four-manifolds with certain properties} \label{examples}

In this section, we establish the existence of symplectic
four-manifolds with certain properties which we need in Section
\ref{construction} for our construction of counter examples. This is
stated precisely in Proposition \ref{ four manifold}, which has
appeared in different guises in
\cite{yan;hodge-structure-symplectic}
  and \cite{R. Gompf} and depends on an idea of Johnson and Rees
\cite{JR}.

  \begin{definition}\label{skew structure} Let $G$ be a discrete group.
  A non-degenerate skew structure on $G$ is a non-degenerate skew
bilinear form
$$ <\,,\,> : H^1(G,\R)\times H^1(G,\R) \rightarrow \R $$ which factors
through the cup product, that is, there exists a linear functional
$\sigma: H^2(G,\R) \rightarrow \R$ so that $\langle a,b
\rangle=\sigma(a \cup b),$
  for all $a,b \in H^1(G,\R)$.
\end{definition}
A finitely presentable group $G$ is called a K\"ahler group if it is
the fundamental group of a closed K\"ahler manifold; otherwise it is
a non-K\"ahler group. It was proved in \cite{JR} that any K\"ahler
group and any of its finite index subgroups must admit a
non-degenerate skew structure.

\begin{lemma}\label{integral class} Let $(N,\omega)$ be a closed, symplectic
$4$-manifold so that $\pi_1(N)$ is a finitely presentable group
which admits a non-degenerate skew structure. Then there exists an
integral class $c$ such that the map $L_c :H^1(N) \rightarrow
H^3(N)$ is an isomorphism.

\end{lemma}

\begin{proof} By elementary homotopy theory there is a natural map $f: N
\rightarrow K(G,1)$ such that the induced homomorphism
$$f^{\ast} : H^{\ast}(G,\R) \rightarrow H^{\ast}(N,\R)$$
is an isomorphism in dimension $1$ and injective in dimension $2$.
Let $\langle\, ,\,\rangle$ be a non-degenerate skew structure on $G$
and $\sigma$ be the corresponding functional on $H^2(G,\R)$. Since
$H^2(G,\R)$ is a subspace of $H^2(N,\R)$, $\sigma$ extends to a
functional $\widetilde{\sigma}$ on $H^2(N,\R)$. By
Poincar$\acute{e}$ duality, there exists a class $c$ such that
$$\tilde{\sigma} (a)=\left(a\wedge c,[N]\right),$$ where $a \in
H^2(N,\R)$ and $[N]$ is the fundamental class of $N$. Suppose $x \in
H^1(N,\R)$ such that $L_c(x)=x \wedge c =0 \in H^3(N,\R)$. Then for
any $y \in H^1(N,\R)$ we have $\tilde{\sigma}(y\wedge x) =
\left((y\wedge x) \wedge c,[N] \right)= \left(y \wedge (x \wedge
c),[N] \right)=0$. Note $\tilde{\sigma}(y\wedge x)= \sigma(y \wedge
x)=\langle y,x \rangle$ we conclude that $\langle y,x \rangle=0$ for
any $y \in H^1(N,R)$. It then follows from the non-degeneracy of
$\langle\,\, ,\,\,\rangle$ that $x=0$.  This shows that $L_c$ is
injective. Then by Poincar$\acute{e}$ duality $L_c$ must be an
isomorphism indeed. Finally note that the set
$$\{\alpha \in \H^2(N) \mid L_{\alpha} :
H^{1}(N) \rightarrow H^{3}(N)\, \text{is an isomorphism} \}$$ is an
open subset of $H^2(N)$. Without the loss of generality, we may
assume that the class $c$ we obtained above is rational. Replace $c$
by $nc$ for some sufficiently large integer $n$ if necessary, we get
an integral class $c$ such that the map $L_c :H^1(N) \rightarrow
H^3(N)$ is an isomorphism. \end{proof}\qed

Combining Lemma \ref{integral class} and Theorem \ref{Gompf}, we get
the existence of symplectic four-manifolds with the desired
properties as stated in the following proposition.

\begin{proposition} \label{ four manifold} Let $G$ be a finitely presentable
group which admits a non-degenerate skew structure. Then there is a
closed, symplectic $4$-manifold $(N,\omega)$ with $\pi _1(N)=G$ such
that the following two conditions are satisfied:
\begin{itemize}
\item[1] the Lefschetz map $L_{[\omega]} : H^1(N) \rightarrow
H^3(N)$ is identically zero. \item[2] there exists an integral class
$c \in H^2(N)$ such that the map $L_c :H^1(N) \rightarrow H^3(N)$ is
an isomorphism. \end{itemize}
\end{proposition}

\section{A remark on the fundamental groups of  strong
Lefschetz four-manifolds}\label{fundamental group}

As an application of Proposition \ref{ four manifold}, we record in
this section an interesting observation on the fundamental groups of
strong Lefschetz four-manifolds.

Using the Hard Lefschetz theorem, Johnson and Rees proved in
\cite{JR} that if a finitely presentable group $G$ is the
fundamental group of a compact K\"ahler manifold, then $G$ has to
admit a non-degenerate skew structure. We note that the fundamental
groups of strong Lefschetz manifolds also have to admit a
non-degenerate skew structure, and Johnson and Rees's argument
applies verbatim to our situation. On the other hand, if $G$ is a
finitely presentable group which supports a non-degenerate skew
structure, then by Proposition \ref{ four manifold} there exists a
compact symplectic four manifold $(N, \omega_0)$ and a closed two
form $c$ on $N$ such that the Lefschetz map $L_{[\omega_0]} : H^1(N)
\rightarrow H^3(N)$ is identically zero and such that the map $L_c
:H^1(N) \rightarrow H^3(N)$ is an isomorphism. For a sufficiently
small constant $\epsilon>0$, set $\omega'=\omega_0+\epsilon c$. It
is easy to see that $\omega'$ is symplectic and satisfies the strong
Lefschetz property. In summary we have the following result.

\begin{theorem} \label{characterization}Suppose $G$ is a finitely
presentable group. Then the following statements are equivalent:
\begin{enumerate}\item $G$ admits a non-degenerate skew structure.
\item $G$ can be realized as the fundamental group of a compact
strong Lefschetz four manifold .
\end{enumerate}
\end{theorem}

Theorem \ref{characterization} raises a natural question whether
there exist finitely presentable non-K\"ahler groups which support a
non-degenerate skew structure. This question is answered
affirmatively in Lemma \ref{exotic group}. However, to prove this
lemma we will need a non-trivial fact concerning non-K\"ahler groups
which is due to Johnson and Rees.

\begin{theorem}\label{grouptheoretic}\cite{JR} Let $G_1$, $G_2$ be groups
which both have at least one nontrivial finite quotient, and let $H$
be any group. Assume that $G=(G_1 \ast G_2)\times H$ admits a
non-degenerate skew structure. Then $G$ has a subgroup of finite
index which does not support any non-degenerate skew structure, and
consequently is not a K\"ahler group itself.
\end{theorem}

\begin{lemma} \label{exotic group}For any positive composite number $m,n$,
the group $G_{m,n}= (\Z_m \ast \Z_n)\times (\Z \times \Z)$ admits a
non-degenerate skew structure; furthermore $G_{m,n}$ has a subgroup
of finite index which does not admit any non-degenerate skew
structure, and therefore is not a K\"ahler group.
\end{lemma}

\begin{proof} Since $m, n$ are composite numbers, both $\Z_m$ and $\Z_n$
have nontrivial finite quotient. It follows from Theorem
\ref{grouptheoretic} that the group $G_{m,n}$ has a subgroup of
finite index which does not support any non-degenerate skew
structure.  Note that by corollary $6.2.10$ and exercise $6.2.5$ of
\cite{Charles A. Weibel}, $H^i(\Z_m\ast \Z_n, \R)=H^i(\Z_n,
\R)\oplus H^i(\Z_m, \R)=0$ for $i \geqslant 1$. Then it follows from
the K$\ddot{u}$nneth formula in group cohomology( see for instance
exercise $6.1.10$ of \cite{Charles A. Weibel}) that $H^{i}(G_{m,n},
\R)=H^i(\Z\times \Z, \R)$ for $i \geqslant 1$. Since $(\Z \times
\Z)$ is a K\"ahler group, $(\Z \times \Z)$ must have a
non-degenerate skew structure. It follows that $G_{m,n}$ also has
such a structure.  \end{proof}\qed

\begin{example} \label{covering} Let $m, n$ be two composite natural numbers
and let $G_{m,n}$ be defined as in Lemma \ref{exotic group}. Since
$G_{m,n}$ does support a non-degenerate skew structure itself, by
Theorem \ref{characterization} it can be realized as the fundamental
group of some symplectic four manifold $N$. By Lemma \ref{exotic
group} $G_{m,n}$ must have a subgroup $K$ of finite index which does
not support any non-degenerate skew structure at all. Let
$\widetilde{N}$ be the finite covering space of $N$ with fundamental
group $K$. Then by Theorem \ref{characterization} again we have that
$\widetilde{N}$ does not support any symplectic form $\omega$ such
that $(\widetilde{N}, \omega)$ has the strong Lefschetz property.
\end{example}

Gompf proved in \cite{R. Gompf} the remarkable result that any
finitely presentable group can be realized as the fundamental group
of a symplectic four-manifold. In contrast, Theorem
\ref{characterization} imposes a rather stringent restriction on the
fundamental groups of compact strong Lefschetz four manifolds. For
example, any non-trivial finitely presentable free group can not be
the fundamental group of a compact strong Lefschetz four manifold.
(C.f., page 592-593 of \cite{R. Gompf}.) In addition, Theorem
\ref{characterization} also asserts that, different from the
fundamental groups of compact K\"ahler manifolds to which far more
rich restrictions apply (see e.g., \cite{AB96}), the fundamental
groups of  compact strong Lefschetz four manifolds have only one
restriction as we stated in Theorem \ref{characterization}.
Therefore, as suggested by Example \ref{covering}, fundamental
groups may serve as effective tools to distinguish strong Lefschetz
manifolds from K\"ahler manifolds.

\section{Examples that the strong Lefschetz property is not
preserved by symplectic reduction}  \label{construction}

Since in this section we are going to make an extensive use of  the
Leray-Hirsch theorem, we first give its precise statement here and
refer to \cite{bott-tu;differential-forms} for details.

\begin{theorem}[Leray-Hirsch theorem]\label{Leray-Hirsch}
Let $E$ be a fiber bundle over $M$ with fiber $F$. Suppose $M$ has a
finite good cover\footnote{An open cover $\{U_{\alpha}\}_{\alpha \in
\Lambda}$ of an $n$ dimensional manifold $M$ is called a good cover
if all non-empty finite intersection $U_{\alpha_0} \cap \cdots \cap
U_{\alpha_p}$ is diffemorphic to $\R^n$. It is well-known that every
compact manifold has a finite good cover. See e.g.,
\cite{bott-tu;differential-forms}}. If there are global cohomology
classes $e_1, e_2, \cdots, e_r$ which when restricted to each fiber
freely generate the cohomology of the fiber, then $H^{\ast}(E)$ is a
free module over $H^{\ast}(M)$ with basis $\{e_1, e_2, \cdots,
e_r\}$, i.e.,
$$H^{\ast}(E) \backsimeq H^{\ast}(M)\otimes \R\{e_1,e_2, \cdots,
e_r\} \backsimeq H^{\ast}(M)\otimes H^{\ast}(F).$$
\end{theorem}

The following proposition enables us to construct six-dimensional
Hamiltonian symplectic manifolds which have the strong Lefschetz
property from the symplectic four-manifolds with properties stated
in Proposition \ref{ four manifold}.

\begin{proposition} \label{spherebundle} Suppose $(N,\omega_0)$ is a
$4$-dimensional compact symplectic manifold  such that:
\begin{enumerate}
\item  the Lefschetz map $L_{[\omega_0]} : H^1(N) \rightarrow
H^3(N)$ is not an isomorphism. \item there exists an integral
cohomology class $[c] \in H^2(N)$ such that the map $L_{[c]} :H^1(N)
\rightarrow H^3(N)$ is an isomorphism.
\end{enumerate}
Then there exists a $S^2$ bundle $\pi_M :M \rightarrow N$ which
satisfies the following conditions:
\begin{enumerate}

\item there is a symplectic form $\omega$ on $M$ such that $(M, \omega)$ has
the strong Lefschetz property;

\item  there is an $S^1$ action on $M$  such that $(M,\omega,S^1)$ is a
compact Hamiltonian manifold which has a non-Lefschetz symplectic
quotient.

\end{enumerate}

\end{proposition}

\begin{proof}Let $S^2$ be the set of unit vectors in $R^3$. In
cylindrical polar coordinates $(\theta, h)$ away from the poles,
where $0\leq \theta <2\pi, -1\leq h \leq 1$,  the standard
symplectic form on $S^2$ is the area form given by $\sigma= \theta
\wedge d h$ . The circle $S^1$ acts on $(S^2, \sigma)$ by rotations
$$ e^{it}(\theta, h)= (\theta+t, h).$$  This action is Hamiltonian
with the moment map given by $\mu=h$, i.e., the height function.

Let $\pi_P: P \rightarrow N $ be the principle $S^1$ bundle with
Euler class $[c]$, let $\Theta$ be the connection $1$-form such that
$d\Theta =\pi^*_P c$, and let $M$ be the associated bundle
$P\times_{S^1} S^2$. Then $\pi_M: M \rightarrow N$ is a symplectic
fibration over the compact symplectic four-manifold $N$. The
standard symplectic form $\sigma$ on $S^2$ gives rise to a
symplectic form $\sigma_x$ on each fibre $\pi_M^{-1}(x)$, where
$x\in N$.  The $S^1$-action on $S^2$ that we described above induces
a fibrewise $S^1$ action on $M$. Furthermore, there is a globally
defined function $H$ on $M$ such that the restriction of $H$ to each
fiber $S^2$ is just the height function $h$.

Next, we resort to minimal coupling construction to get a closed two
form $\eta$ on $M$ which restricts to the forms $\sigma_x$ on the
fibres. Let us give a sketch of this construction here and refer to
\cite{A. Weinstein}, \cite{R. Sternberg} and \cite{GS1} for
technical details. Consider the closed two form $-d(t \Theta)=
-td\Theta-dt \wedge \Theta$ defined on $P \times \R$. It is easy to
see the $S^1$ action on $P \times \R$ given by
$$e^{it}(p, t)=(e^{it}p, t)$$ is Hamiltonian with the moment map
$t$. Thus the diagonal action of $S^1$ on $(P \times \R) \times S^2$
is also Hamiltonian, and $M$ is just the reduced space of $(P\times
\R)\times S^2$ at the zero level. Moreover, the closed two form
$\left(-d(t\Theta) +\sigma \right) \mid _{\text{zero level}}$
descends to a closed two form $\eta$ on $M$ with the desired
property.

It is useful to have the following explicit description of $\eta$.
Observe that $\theta-\Theta$ is a basic form on $(P \times \R)
\times S^2$. Its restriction to the zero level of $(P \times \R)
\times S^2$ descends to a one form $\tilde{\theta}$ on $M$ whose
restriction to each fibre $S^2$ is just $\theta$. It is easy to see
that on the associated bundle $P\times_{S^1}(S^2-\text{\{two
poles\}})$ we actually have $\eta = H\pi^*_M c+ dH \wedge
\tilde{\theta}$.

For any real number $t_0$, note that the restriction of $\eta -t_0
\pi_M ^* c$  to fibres are symplectic forms $\sigma_x$. By an
argument due to Thurston \cite{MS98},  for sufficiently large
constant $K>0$ the form $K\pi_M^*\omega_0 -t_0 \pi_M^* c+\eta$ is
symplectic. Equivalently, define $ \omega=\pi_M^*\omega_0 -\epsilon
t_0 \pi_M^*c+ \epsilon \eta$ for sufficiently small constant
$\epsilon>0$. Then $\omega$ is a symplectic form on $M$;
furthermore, the fibrewise $S^1$-action on $(M, \omega)$ is
Hamiltonian with the moment map $H: M \rightarrow \R$.

Choose some $min_{x \in M} H(x)<t_0< max_{x\in M}H(x)$ and have it
fixed. If we perform symplectic reduction at $H =t_0$,  the
symplectic reduced space is $N$ with the reduced form $\omega_0$.
Clearly, $(N, \omega_0)$ does not satisfy the strong Lefschetz
property since the Lefschetz map $L_{[\omega_0]}: H^1(N) \rightarrow
H^3(N)$ is not an isomorphism.

It remains to check that for sufficiently small constant
$\epsilon>0$, $(M, \omega)$ has the strong Lefschetz property.

Consider the closed $2$-form $\eta$ on $M$. Its restriction to each
fibre $S^2$ generates the second cohomology group $H^2(S^2)$. Write
$H(S^2)=\R[x]/(x^2)$, where $\R[x]$ is the real polynomial ring and
$(x^2)$ is the ideal of $\R[x]$ generated by the quadratic
polynomial $x^2$.  By the Leray-Hirsch theorem there is an additive
isomorphism
$$H(N)\otimes \R[x]/(x^2) \rightarrow H(M), \, [\alpha] \otimes x^i
\rightarrow [\pi_M^{\ast}\alpha \wedge \eta^i],\, i=0,1.$$ As a
result we have $[\eta ^2]=[\pi_M^{\ast}\beta_2 \wedge \eta]+
[\pi_M^{\ast}\beta_4]$, where $\beta_2$ and $\beta_4$ are closed
forms on $N$ of degree two and four respectively.

Choose an $\epsilon>0$ which is sufficiently small such that
\begin{equation}\label{neq1}
[\omega_{0}-t_0 \epsilon c]^2 \neq -\epsilon^2 [\beta_4] + \epsilon
[(\omega_0 -t_0 \epsilon c) \wedge \beta_2].
\end{equation}
We claim for the $\epsilon$ chosen above, the symplectic manifold
$(M, \pi_M^*\omega_0 -\epsilon t_0 \pi_M^*c+ \epsilon \eta)$ will
satisfy the strong Lefschetz property. By Poincar$\acute{e}$ duality
it suffices to show the two Lefschetz maps
\begin{equation}\label{injective map1}
L^2_{[\omega ]}: H^1(M) \rightarrow H^5(M)
\end{equation}
\begin{equation}\label{injective map2}
L_{[\omega ]}: H^2(M) \rightarrow H^4(M)
\end{equation}
are injective. We will give a proof in two steps below.
\begin{enumerate}

\item It follows from the Leray-Hirsch theorem that $H^1(N)
\xrightarrow [\pi_M^{\ast}]{\backsimeq} H^1(M)$. Thus to show Map
(\ref{injective map1}) is injective we need only to show for any $
[\lambda] \in H^1(N)$ if $L^2_{[\omega ]} (\pi_M^{\ast}[\lambda])=0$
then we have $[\lambda]=0$. Since $\omega= \pi_M^{\ast}(\omega_0
-t_0\epsilon c) +\epsilon \eta$, $[\eta^2] =[\pi_M^{\ast}\beta_2
\wedge \eta]+ [\pi_M^{\ast}\beta_4]$ and any forms on $N$ with
degree greater than $4$ vanishes, we have
\begin{equation}
\begin{split} 0& = L^2_{[\omega ]} ([\pi^{\ast}_M\lambda])\\
               &= \pi_M^{\ast}\left(2\epsilon
[\omega_0-t_0\epsilon c] + \epsilon^2 [\beta_2]\right)
\wedge[\pi_M^{\ast} \lambda]\wedge [\eta].\end{split}\end{equation}
Since by the Leray-Hirsch theorem $H(M)$ is free over $1$ and
$[\eta]$, we get that \[\begin{split}  0 &= \pi_M^{\ast}
(2\epsilon[\omega_0]-2t_0\epsilon^2 [c] + \epsilon^2 [\beta_2])
\wedge \pi^{\ast}[\lambda]\\&= \pi_M^{\ast}\left(
([2\epsilon\omega_0-2t_0\epsilon^2 c + \epsilon^2 \beta_2]) \wedge
[\lambda]\right).\end{split}\]   Since $L_{[c]}: H^1(N) \rightarrow
H^3(N)$ is an isomorphism, the determinant of the linear map
$L_{[2\epsilon\omega_0-2t_0\epsilon^2 c + \epsilon^2 \beta_2]}:
H^1(N) \rightarrow H^3(N)$ is a polynomial in $t_0$ of positive
degree. Therefore $L_{[2\epsilon\omega_0-2t_0\epsilon^2 c +
\epsilon^2 \beta_2]}: H^1(N) \rightarrow H^3(N)$ is an isomorphism
except for finitely many possible values of $t_0$. If necessary,
replace $H$ and $t_0$ by $H+c$ and $t_0+c$ respectively for some
suitable small constant $c>0$. We conclude that Map (\ref{injective
map1}) is an isomorphism.

\item By the Leray-Hirsch theorem, to show that Map
(\ref{injective map2}) is injective it suffices to show if
$L_{[\omega]}(\pi_M^{\ast}[\varphi]+ k [\eta])=0$ for arbitrarily
chosen scalar $k$ and second cohomology class $[\varphi] \in
H^2(N)$, then we have $[\varphi]=0$ and $k=0$. Since
$\omega=\pi_M^{\ast}(\omega_0-t_0\epsilon c) +\epsilon \eta$ and
$[\eta^2] =[\pi_M^{\ast}\beta_2 \wedge \eta]+
[\pi_M^{\ast}\beta_4]$, we have
\begin{equation} \begin{split} 0&=L_{[\omega]}(\pi_M^{\ast}[\varphi]+ k
[\eta])\\&=\left( \pi_M^{\ast}[(\omega_0-t_0\epsilon
c)\wedge\varphi] + \epsilon k \pi_M^{\ast}[\beta_4]
\right)+\\&\,\,\,\,\,\, \left( k \pi_M^{\ast}[\omega_0-t_0\epsilon c
] + \epsilon \pi_M^{\ast}[\varphi] +\epsilon k \pi_M^{\ast}[\beta_2]
\right)\wedge \eta \end{split}
\end{equation} By the Leray-Hirsch theorem $H(M)$ is a free
module over $1$ and $[\eta]$. So we have that
\begin{equation} \label{part1}
\pi_M^{\ast}[(\omega_0-t_0\epsilon c)\wedge\varphi] + \epsilon k
\pi_M^{\ast}[\beta_4]=0
\end{equation}
\begin{equation}\label{part2}
k \pi_M^{\ast}[\omega_0-t_0\epsilon c ] + \epsilon
\pi_M^{\ast}[\varphi] +\epsilon k \pi_M^{\ast}[\beta_2]=0
\end{equation}
If $k=0$, it follows easily from the equation (\ref{part2}) that
$[\varphi]=0$. Assume $k \neq 0$. Substitute $
\pi_M^{\ast}[\varphi]= - \dfrac{1}{\epsilon} k
\pi_M^{\ast}[\omega_0-t_0\epsilon c]-k \pi_M^{\ast}[\beta_2]$ into
the equation (\ref{part1}) we get
$$ \pi_M^{\ast}[\omega_0-t_0\epsilon
c]\wedge (-k \pi_M^{\ast}[\omega_0-t_0\epsilon c]-\epsilon k
\pi_M^{\ast}[\beta_2])+ \epsilon^2 k \pi_M^{\ast}[\beta_4]=0
$$
Since $k\neq 0$, we get
$$\pi_M^{\ast}([\omega_0-t_0\epsilon c])^2 = -\epsilon^2\pi_M^{\ast}
[\beta_4] + \epsilon \pi_M^{\ast}[(\omega_0-t_0\epsilon c) \wedge
\beta_2]$$ This contradicts the equation (\ref{neq1}).
\end{enumerate}
\end{proof}\qed

Now we are in a position to construct examples that the strong
Lefschetz property does not survive symplectic reduction.

\begin{example} \label{non-invariant-kahler}Since the torus is a K\"ahler
manifold, $G=\Z\times \Z$ is a K\"ahler group and thus admits a
non-degenerate skew structure. Clearly, by Lemma \ref{ four
manifold} there is a closed, symplectic $4$-manifold $(N,\omega_0)$
which satisfies the following conditions:
\begin{enumerate} \item $\pi _1(N)=\Z\times \Z$ \item The
Lefschetz map $L_{[\omega]} :H^1(N) \rightarrow H^3(N)$ is trivial.
\item There is an integral class $[c]\in H^2(N)$ such that the map
$L_{[c]}:H^1(N) \rightarrow H^3(N)$ is an isomorphism.
\end{enumerate} Then it follows easily from Proposition
\ref{spherebundle} that there exists a compact six-dimensional
Hamiltonian circle manifold $(M, \omega)$ which has the strong
Lefshetz property itself but admits a non-Lefschetz symplectic
quotient.
\end{example}

Since the six-dimensional Hamiltonian $S^1$-manifold $(M, \omega)$
constructed in Example \ref{non-invariant-kahler} has a
non-Lefschetz symplectic quotient, $\omega$ can not be an invariant
K\"ahler form. But in general we do not know whether $M$ supports
any K\"ahler form or not. To get examples which do not admit any
K\"ahler structure, we observe that by our construction $M
\rightarrow N$ is a fibration with fiber $S^2$ and so $\pi_1(M) =
\pi_1(N)$. Instead of choosing $G=\Z \times \Z$, we may well choose
$G=G_{m,n}$, where $m, n$ are any composite numbers and $G_{m,n}$ is
defined as in Lemma \ref{exotic group}. For any such a group $G$,
the corresponding Hamiltonian manifold $M$ has a non-K\"ahler
fundamental group and therefore is not homotopy equivalent to any
compact K\"ahler manifold. Thus we have proved the following
theorem:

\begin{theorem} \label{counter} There exist infinitely many
topologically inequivalent six-dimensional compact Hamiltonian
symplectic $S^1$-manifolds which satisfy the following conditions:

\begin{enumerate} \item
the strong Lefschetz property,

\item admitting a non-Lefschetz symplectic quotient, \item  not
homotopy equivalent to any compact K\"ahler manifold.
\end{enumerate}
\end{theorem}

Finally we observe that the fixed point set of the Hamiltonian
symplectic manifold $(M, \omega, S^1)$ constructed in Proposition
\ref{spherebundle} has two components on which the moment map takes
maximum and minimum respectively; furthermore, in the proof of
Proposition \ref{spherebundle}, if we choose $t_0$ to be the minimum
value of the moment map, then the minimal component as a symplectic
submanifold can be identified with $(N, \omega_0)$ which clearly
does not have the strong Lefschetz property. This observation,
together with Lemma \ref{exotic group}, leads to the following
result.

\begin{theorem}There exist infinitely many topologically
inequivalent six-dimensional compact Hamiltonian symplectic
$S^1$-manifolds which satisfy the following conditions:

\begin{enumerate} \item
the strong Lefschetz property,

\item admitting a non-Lefschetz fixed point sub-manifold, \item  not
homotopy equivalent to any compact K\"ahler manifold.
\end{enumerate}

\end{theorem}




\end{document}